\documentclass[12pt,a4paper]{amsart}
\usepackage{amssymb,amsfonts,amsthm,amsmath,graphicx,color}

\theoremstyle{plain}

\numberwithin{equation}{section} \numberwithin{theorem}{section}
\numberwithin{lemma}{section} \numberwithin{definition}{section}
\numberwithin{corollary}{section}
\numberwithin{proposition}{section} \textheight =24cm
\begin{document}
\title[The series that Ramanujan misunderstood]{The series that Ramanujan misunderstood}
\author{Geoffrey B Campbell}
\address{Mathematical Sciences Institute,
         The Australian National University,
         Canberra, ACT, 0200, Australia}

\email{Geoffrey.Campbell@anu.edu.au}

\author{Aleksander Zujev}
\address{Department of Physics,
         University of California,
         Davis, CA, USA}

\email{azujev@ucdavis.edu}


\begin{abstract}
We give a new appraisal of the function $\Delta(x)$ and its zeroes in the equation $f(x) = g(x) + \Delta(x)$ where $f(x) = \sum_{n \in Z} 2^n x^{2^n}$ and $g(x) = 1/((\log2)(\log(1/x)))$.
\end{abstract}

\maketitle

\section{Introduction} \label{S:intro}

Consider the bilateral infinite series that converges in the unit disc $|x|<1$,

\begin{eqnarray} \label{eq1}
&& f(x) = \sum_{n \in Z} 2^n x^{2^n} \\
&&  =... + 1024x^{1024} + 512x^{512} + 256x^{256} + 128x^{128} + 64x^{64} + 32x^{32} + 16x^{16} \nonumber \\
&& + 8x^8 + 4x^4 + 2x^2 + x + \frac{x^{1/2}}{2} + \frac{x^{1/4}}{4} + \frac{x^{1/8}}{8} + \frac{x^{1/16}}{16} \nonumber \\
&&  + \frac{x^{1/32}}{32} + \frac{x^{1/64}}{64} + \frac{x^{1/128}}{128} + \frac{x^{1/256}}{256} \nonumber \\
&&  + \frac{x^{1/512}}{512} + \frac{x^{1/1024}}{1024} + ... \nonumber
\end{eqnarray}
Next also for $|x|<1$, consider the function,
\begin{equation} \label{eq2}
 g(x) = \frac{1}{\log2 \log(1/x)}
\end{equation}

 Ramanujan, in his theory of prime numbers in his pre-Cambridge days, seemed to believe that for all real $0<x<1, f(x) = g(x)$. In Hardy's famous book on Ramanujan \cite{Hardy1940}, we can form a view that Ramanujan was familiar with the Euler-McLaurin summation formula from the Carr Synopsis book he referred to constantly, and that this formula omitted the oscillating term $\Delta(x)$. As a result, Ramanujan inferred many things about the distribution of prime numbers as if there were no analytic theory introduced by Riemann in his landmark paper of 1859 which put the now-named Riemann Hypothesis, and gave the first proof of the Riemann zeta functional equation. Using contour integration and the residue theorem, the reality is that
\begin{equation}  \label{eq3}
 f(x) = g(x) + \Delta(x),
\end{equation}

where $\Delta(x)$ oscillates   
around zero, and the amplitude of the oscillations only are noticeable around the third or fourth decimal place.
Indeed, both $f(x)$ and $g(x)$ satisfy the functional equation $2f(x^2) = f(x)$, and $f(x) = g(x)$
approximately to 3 or 4 decimal places. The $\Delta(x)$ oscillations become “more wriggly” as $x$ approaches $1$
 near it’s limiting boundary value of convergence. G H Hardy was able to explain to Ramanujan that
 $\Delta(x)$ is an oscillating periodic function of $\log(\log(1/x))$.
The correct formula corresponding to (\ref{eq3}) is for $|x|<1$,
\begin{equation}  \label{eq4}
 \sum 2^k x^{2^k} = \frac{1}{\log2 \log(1/x)} \left\{1 - {\sum} ' \Gamma\left(1+ \frac{2ki\pi}{\log 2}\right) \left(\log\left(\frac{1}{x}\right)\right)^{-2ki\pi/\log 2} \right\},
\end{equation}
with the sum $\sum$ over all integers $k$, and the sum $\sum ' $ over all nonzero integers $k$.

The problem is to locate the zeroes of $\Delta(x)$, and so find where (\ref{eq3}) above becomes $f(x) = g(x)$.

\section{Approximation with self-similar oscillating function}

\subsection{Function $\Delta_0(x)$.}

\begin{figure}
\includegraphics[width=3in,angle=0]{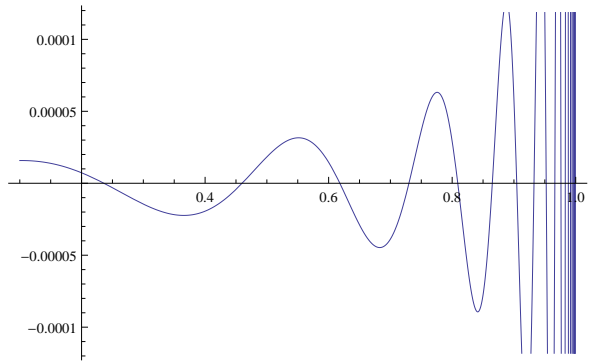}
\caption{(Color online)
Plot of $\Delta_0(x)$.
As $x$ approaches $1$, with every oscillation, frequency of oscillations of $\Delta_0(x)$ increases 2 times,
and its amplitude increases 2 times.
\label{Delta0}}
\end{figure}

At $x$ close to $1$, $\log(1/x) \; \approx \; (1-x)$, and
\begin{eqnarray}
&&  \label{eq5} \Delta(x) = \frac{1}{\log 2 \log(1/x)} {\sum} ' \Gamma\left(1+ \frac{2ki\pi}{\log 2}\right) (\log(1/x))^{-2ki\pi/\log 2} \approx \\
&&  \label{eq6}  \Delta_0(x) = \frac{1}{\log 2} {\sum} ' \Gamma\left(1+ \frac{2ki\pi}{\log 2}\right) (1-x)^{-1-2ki\pi/\log 2}
\end{eqnarray}
$\Delta_0(x)$ is a self-similar function, such that $\Delta_0((x + 1)/2) = 2\Delta_0(x)$.
As $x$ approaches $1$, period of oscillations of $\Delta_0(x)$ exponentially decreases,
and its amplitude exponentially increases.
Fig. (\ref{Delta0}) shows the plot of $\Delta_0(x)$.
Due to such periodicity of $\Delta_0(x)$, it is enough to study this function at any interval $[x, (x+1)/2]$
for complete knowledge of the function.
The function $\Delta_0(x)$ is dominated by the largest ($k = \pm 1$) terms of the sum (\ref{eq6}),
and these two terms add to a function of the form
\begin{equation}
\label{eq7}  \frac{b}{1-x}  \cos \left( \log(1-x) \frac{2\pi}{\log 2} +\phi \right) .
\end{equation}
It is sinusoide, which get squeezed horizontally as $x$ approaches $1$, and get stretched vertically.

We can study and write more about the function $\Delta_0(x)$ if needed. In the paper by Campbell \cite{Campbell1994} he refers to an ingenious approach to finding zeroes of a similar oscillating function examined in a study by Mahler \cite{Mahler1980}, which may be applicable for the functions in our current paper.

Of particular interest to us are zeroes of $\Delta_0(x)$.
The first zero of $\Delta_0(x)$ is
$x_0 \approx 0.23628629$.
All consecutive zeroes are given by
\begin{equation}  \label{Delta0zeron}
 x_n = 1 - \frac{1-x_0}{2^{n/2}}.
\end{equation}

\subsection{Approximation of $\Delta(x)$ by $\Delta_0(x)$.}

\begin{figure}
\includegraphics[width=2.5in,angle=0]{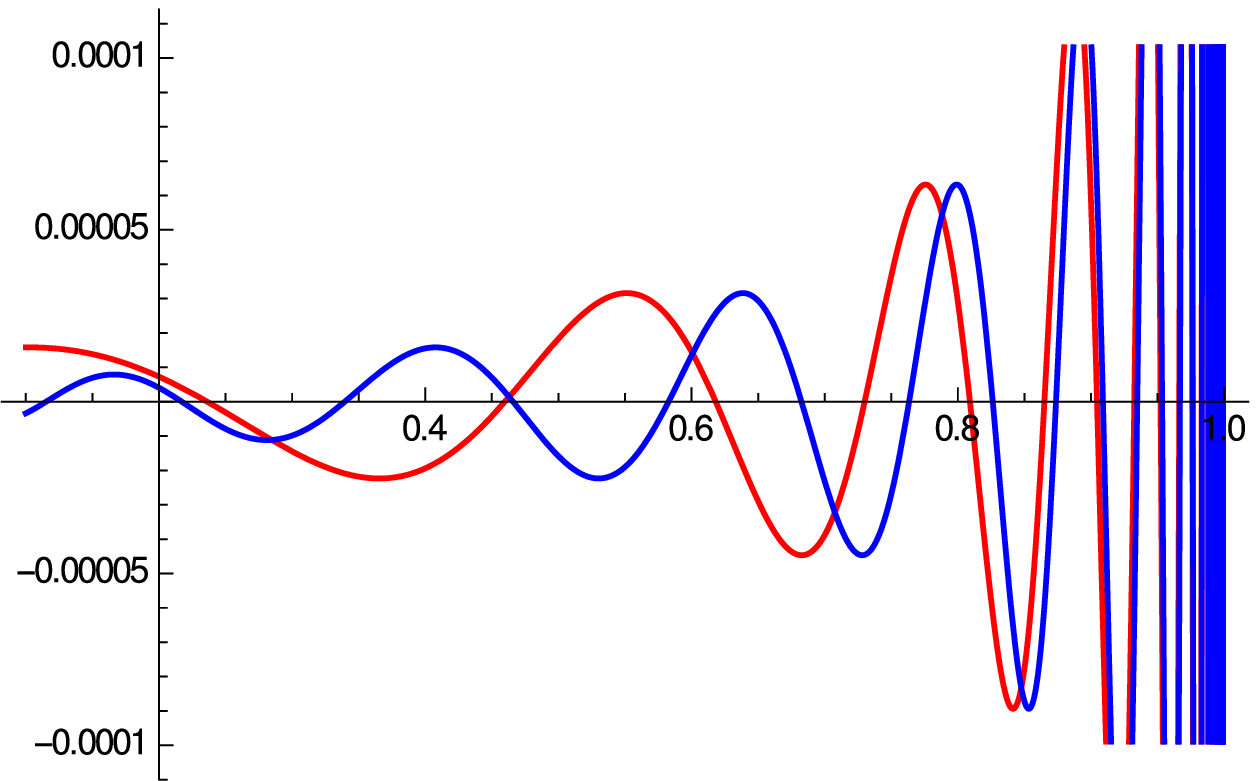}
\includegraphics[width=3in,angle=0]{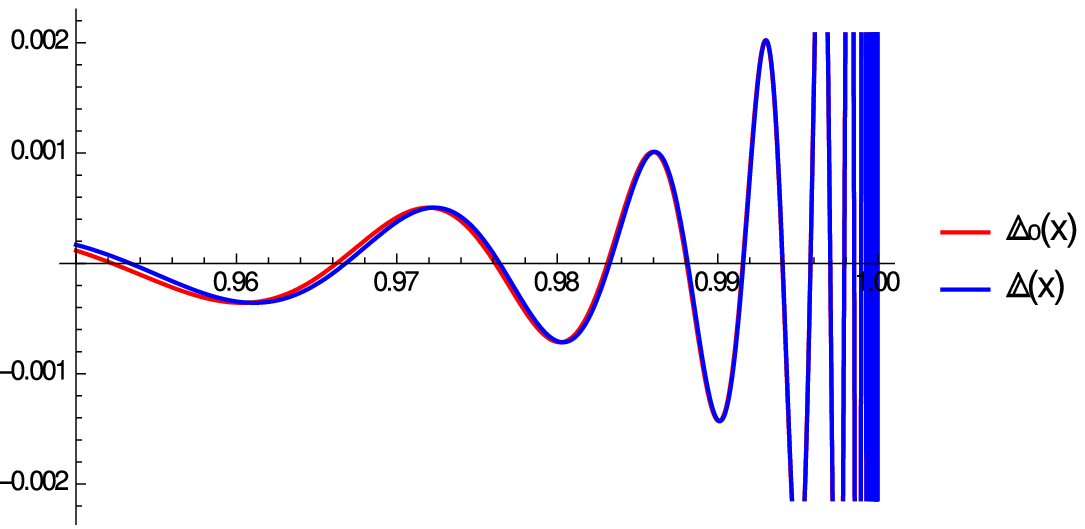}
\caption{(Color online)
Plots of $\Delta_0(x)$ and $\Delta(x)$ at intervals [0.2, 1] (left) and [0.95, 1] (right).
As $x$ approaches $1$, $\Delta_0(x)$ and $\Delta(x)$ converge.
\label{Delta0Delta}}
\end{figure}

How well is $\Delta(x)$ approximated by $\Delta_0(x)$?
Fig. (\ref{Delta0Delta}) shows the plots of both $\Delta_0(x)$ and $\Delta(x)$.
At smaller $x$, $\Delta_0(x)$ and $\Delta(x)$ differ considerably, but as $x$ approaches $1$,
$\Delta_0(x)$ and $\Delta(x)$ converge.

\subsubsection{Numerical estimates.}

A few first zeroes of $\Delta(x)$ and $\Delta_0(x)$,
and relative error of approximation, given by
$|$((zero of $\Delta(x)$) - (zero of $\Delta_0(x)$))/(1 - (zero of $\Delta(x)$))$|$
is shown in Table (\ref{TableDelta}).

\begin{table}
\begin{tabular}{| l | l | l |}
\hline
zero of $\Delta(x)$ & zero of $\Delta_0(x)$ & relative error \\
\hline
0.4659328665 & 0.2362862900 & 0.4299957 \\
0.5827324804 & 0.4599728568 & 0.2941988 \\
0.6825927537 & 0.6181431450 & 0.2030502 \\
0.7633691635 & 0.7299864284 & 0.1410752 \\
0.8261917175 & 0.8090715725 & 0.0985002 \\
0.8737099995 & 0.8649932142 & 0.0690220 \\
0.9089508885 & 0.9045357863 & 0.0484914 \\
0.9347245581 & 0.9324966071 & 0.0341315 \\
0.9533891590 & 0.9522678931 & 0.0240559 \\
0.9668115422 & 0.9662483035 & 0.0169709 \\
0.9764164885 & 0.9761339466 & 0.0119805 \\
0.9832657536 & 0.9831241518 & 0.0084618 \\
0.9881378894 & 0.9880669733 & 0.0059784 \\
0.9915975764 & 0.9915620759 & 0.0042250 \\
0.9940512509 & 0.9940334866 & 0.0029862 \\
0.9957899259 & 0.9957810379 & 0.0021111 \\
0.9970211888 & 0.9970167433 & 0.0014924 \\
0.9978927427 & 0.9978905190 & 0.0010552 \\
0.9985094836 & 0.9985083717 & 0.0007460 \\
0.9989458157 & 0.9989452595 & 0.0005276 \\
0.9992544639 & 0.9992541858 & 0.0003730 \\
0.9994727689 & 0.9994726297 & 0.0002639 \\
0.9996271624 & 0.9996270929 & 0.0001865 \\
0.9997363497 & 0.9997363149 & 0.0001319 \\
0.9998135638 & 0.9998135465 & 0.0000930 \\
0.9998681661 & 0.9998681574 & 0.0000663 \\
0.9999067776 & 0.9999067732 & 0.0000469 \\
0.9999340809 & 0.9999340787 & 0.0000334 \\
0.9999533877 & 0.9999533866 & 0.0000236 \\
0.9999670399 & 0.9999670394 & 0.0000154 \\
0.9999766936 & 0.9999766933 & 0.0000120 \\
0.9999835198 & 0.9999835197 & 0.0000071 \\
0.9999883467 & 0.9999883467 & 0.0000018 \\
\hline
\end{tabular}
\caption{Zeroes of $\Delta(x)$ and $\Delta_0(x)$.
\label{TableDelta}}
\end{table}

As $x$ approaches $1$, the relative error goes to zero.
The estimate of relative error can be given comparing Taylor series expansion for $\Delta(x)$ and $\Delta_0(x)$.
If $x_z$ is a zero of $\Delta(x)$, and $x_{z0}$ is corresponding zero of $\Delta_0(x)$, then
\begin{equation}
 (1-x_{z0}) \approx (1-x_z) + \frac{1}{2}(1-x_z)^2 + \frac{1}{3}(1-x_z)^3,
\end{equation}
or
\begin{equation}
 (1-x_{z}) \approx (1-x_{z0}) - \frac{1}{2}(1-x_{z0})^2 - \frac{1}{3}(1-x_{z0})^3.
\end{equation}


\section{Arbitrary $a$}

The results of the previous section are applicable to the equation with an arbitrary $a$ instead of $2$.
Consider the bilateral infinite series that converges for real $0 < x < 1$,

\begin{eqnarray} \label{eq1gen}
&& f(x) = \sum_{n \in Z} a^n x^{a^n} \\
&&  =... + a^{10}x^{a^{10}} + a^9x^{a^9} + a^8x^{a^8} + a^7 x^{a^7} + a^6 x^{a^6} + a^5x^{a^5} + a^4x^{a^4} \nonumber \\
&& + a^3 x^{a^3} + a^2 x^{a^2} + a x^a + x + \frac{x^{1/a}}{a} + \frac{x^{1/a^2}}{a^2} + \frac{x^{1/a^3}}{a^3} + \frac{x^{1/a^4}}{a^4} \nonumber \\
&&  + \frac{x^{1/a^5}}{a^5} + \frac{x^{1/a^6}}{a^6} + \frac{x^{1/a^7}}{a^7} + \frac{x^{1/a^8}}{a^8} \nonumber \\
&&  + \frac{x^{1/a^9}}{a^9} + \frac{x^{1/a^{10}}}{a^{10}} + ... \nonumber
\end{eqnarray}
Next consider the function given by real $0 < x < 1$,
\begin{equation} \label{eq2gen}
 g(x) = \frac{1}{((\log a) \log(1/x))}.
\end{equation}

\begin{equation}  \label{eq3gen}
 f(x) = g(x) + \Delta(x),
\end{equation}

where $\Delta(x)$ oscillates   
around zero, and the amplitude of the oscillations only are noticeable around the third or fourth decimal place.
Both $f(x)$ and $g(x)$ satisfy the functional equation $a f(x^2) = f(x)$, and $f(x) = g(x)$
approximately to 3 or 4 decimal places. The $\Delta(x)$ oscillations become “more wriggly” as $x$ approaches $1$
 near it’s limiting boundary value of convergence. 
 $\Delta(x)$ is an oscillating periodic function of $\log(\log(1/x))$.
The correct formula corresponding to (\ref{eq3gen}) is for $|x|<1$,
\begin{equation}  \label{eq4gen}
 \sum a^k x^{a^k} = \frac{1}{((\log a)(\log(1/x)))} \left\{1 - {\sum} ' \Gamma\left(1+ \frac{2ki\pi}{\log a}\right) \left(\log\left(\frac{1}{x}\right)\right)^{-2ki\pi/\log a} \right\},
\end{equation}
where the sum $\sum$ is over all integers $k$, and the sum $\sum ' $ is over all nonzero integers $k$.

$\Delta(x)$ may be approximated by
\begin{eqnarray}
\label{eq6gen}  \Delta_0(x) = \frac{1}{\log a} {\sum} ' \Gamma\left(1+ \frac{2ki\pi}{\log 2}\right) (1-x)^{-1-2ki\pi/\log a}
\end{eqnarray}

\begin{figure}
\includegraphics[width=2.5in,angle=0]{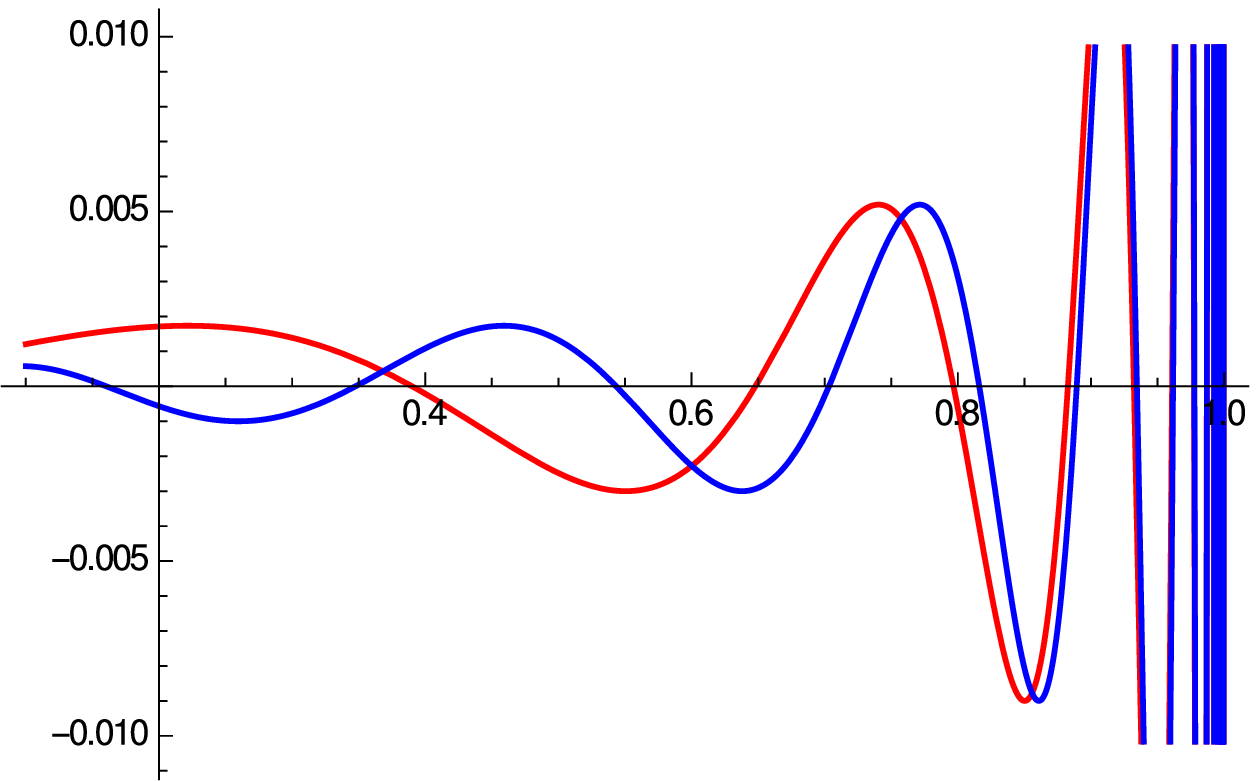}
\includegraphics[width=3in,angle=0]{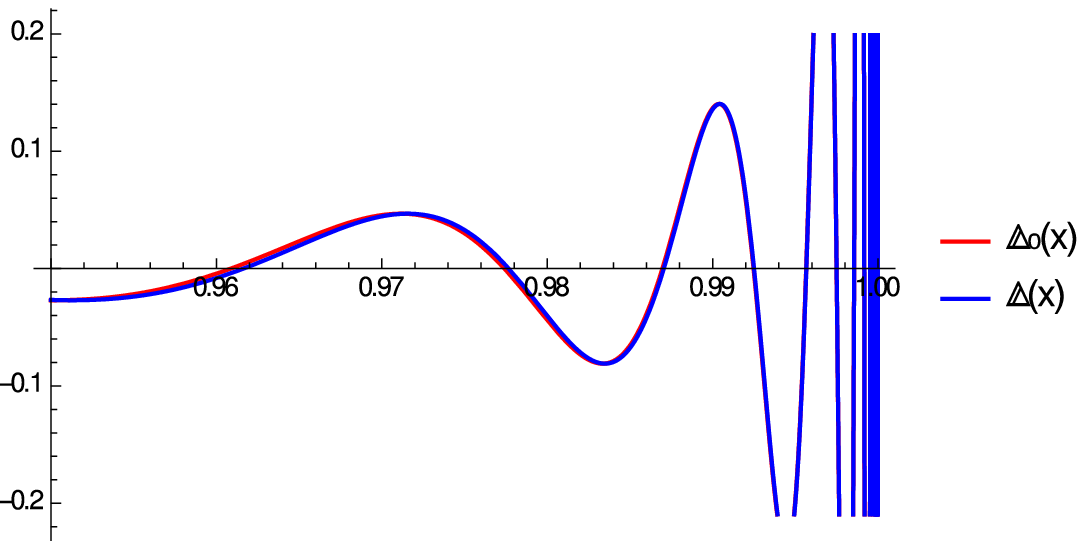}
\caption{(Color online)
Plots of $\Delta_0(x)$ and $\Delta(x)$ at intervals [0.2, 1] (left) and [0.95, 1] (right).
As $x$ approaches $1$, $\Delta_0(x)$ and $\Delta(x)$ converge.
\label{Delta0Deltagen}}
\end{figure}

As an example, we consider $a = 3$.
Fig. (\ref{Delta0Deltagen}) shows the plots of both $\Delta_0(x)$ and $\Delta(x)$.
At smaller $x$, $\Delta_0(x)$ and $\Delta(x)$ differ considerably, but as $x$ approaches $1$,
$\Delta_0(x)$ and $\Delta(x)$ converge.

Zeroes of $\Delta(x)$ may be approximated by zeroes of $\Delta_0(x)$.
The first zero of $\Delta_0(x)$
can be found by numerically solving equation $\Delta_0(x) = 0$.
Approximately, taking only the first terms of the sum,
\begin{equation}
\label{eq7gen}   \Gamma\left(1+ \frac{2i\pi}{\log a}\right) (1-x)^{-1-2i\pi/\log a} + \Gamma\left(1- \frac{2i\pi}{\log a}\right) (1-x)^{-1+2i\pi/\log a} = 0
\end{equation}
\begin{equation}
\label{eq8gen} x_0 \approx 1 - e^{\frac{\left(\frac{\pi}{2} - arg\left(\Gamma\left(1+\frac{2\pi i}{\log a}\right)\right)\right)\log a}{-2\pi}}
\end{equation}

All consecutive zeroes are given by
\begin{equation}  \label{Delta0zerongen}
 x_n = 1 - \frac{1-x_0}{a^{n/2}}.
\end{equation}

\end{document}